# Small-Gain Stability Analysis of Hyperbolic-Parabolic PDE Loops

Iasson Karafyllis[*] and Miroslav Krstic[**]

[*]Dept. of Mathematics, National Technical University of Athens, Zografou Campus, 15780, Athens, Greece, email: iasonkar@central.ntua.gr

[**]Dept. of Mechanical and Aerospace Eng., University of California, San Diego, La Jolla, CA 92093-0411, U.S.A., email: krstic@ucsd.edu

## Abstract

This work provides stability results in the spatial sup norm for hyperbolic-parabolic loops in one spatial dimension. The results are obtained by an application of the small-gain stability analysis. Two particular cases are selected for the study because they contain challenges typical of more general systems (to which the results are easily generalizable but at the expense of less pedagogical clarity and more notational clutter): (i) the feedback interconnection of a parabolic PDE with a first-order zero-speed hyperbolic PDE with boundary disturbances, and (ii) the feedback interconnection, by means of a combination of boundary and in-domain terms, of a parabolic PDE with a first-order hyperbolic PDE. The first case arises in the study of the movement of chemicals underground and includes the wave equation with Kelvin-Voigt damping as a subcase. The second case arises when we apply backstepping to a pair of hyperbolic PDEs that is obtained by ignoring diffusion phenomena. Moreover, the second case arises in the study of parabolic PDEs with distributed delays. In the first case, we provide sufficient conditions for ISS in the spatial sup norm with respect to boundary disturbances. In the second case, we provide (delay-independent) sufficient conditions for exponential stability in the spatial sup norm.

**Keywords:** ISS, parabolic PDEs, hyperbolic PDEs, boundary disturbances.

## 1. Introduction

The use of the notion of Input-to-State Stability (ISS) for finite-dimensional systems, which was proposed by E. D. Sontag in [28], allowed the development of small-gain theorems. Starting with the first nonlinear, generalized small-gain theorem in [11] for systems described by Ordinary Differential Equations (ODEs), the small-gain stability analysis has been extended to various kinds of systems with inputs (see [12]). The extension of ISS to systems described by Partial Differential Equations (PDEs; see [22,23,24]) allowed the application of small-gain arguments in [21] for systems of interconnected PDEs. The recent extension of ISS to PDEs with boundary disturbances in [13,14] allowed the use of small-gain arguments in [14] to PDEs with non-local boundary conditions. The use of small-gain arguments in [14] also showed that small-gain analysis is capable of providing stability estimates in the spatial sup norm. This feature can rarely be met in Lyapunov analysis (which is more well-suited for estimates in $L^p$ spatial norms with $2 \leq p < +\infty$).

The study of interconnected PDEs arises naturally in many applications. Interconnections of PDEs have been studied in [18]. The literature focuses on the study of systems of hyperbolic PDEs (see [3,15,31,32]) and Reaction-Diffusion systems (i.e., systems of parabolic PDEs; see for instance [27]). ODE-PDE cascades have been studied in [1,2,16,17,29], mostly for feedback and observer design purposes. However, the present work is devoted to the study of parabolic-hyperbolic PDE loops. Such loops present unique features because they combine the finite signal transmission speed of hyperbolic PDEs with the unlimited signal transmission speed of parabolic PDEs. Since there are many possible interconnections that can be considered, it is difficult to give results for a "general case". Therefore, we focus on two particular cases, which are analyzed in detail, because they contain challenges typical of more general systems (to which the results are easily generalizable but at the expense of less pedagogical clarity and more notational clutter).

The first case considered in this paper is the feedback interconnection of a parabolic PDE with a special first-order hyperbolic PDE: a zero-speed hyperbolic PDE. Thus the action of the hyperbolic PDE resembles the action of an infinite-dimensional, spatially parameterized ODE. However, the study of this particular loop is of special interest because it arises in an important application: the movement of chemicals underground ([8], pages 210-216). Moreover, the study of this particular system can be used for the analysis of wave equations with Kelvin-Voigt damping (see also [5,7,9,30]). In this case, we provide sufficient conditions for ISS in the spatial sup norm with respect to boundary disturbances (Theorem 2.2 and Corollary 2.3). There are no available stability results in the literature for the wave equation with Kelvin-Voigt damping in the spatial sup norm (even when boundary disturbances are absent).

The second case considered in this paper is the feedback interconnection, by means of a combination of boundary and in-domain terms, of a parabolic PDE with a first-order hyperbolic PDE. The interconnection is effected by linear, non-local terms. The second case arises when we apply backstepping to a pair of hyperbolic PDEs that is obtained by ignoring diffusion phenomena (see [6,10]). Moreover, the second case arises in the study of parabolic PDEs with distributed delays of trace terms. Parabolic equations with delayed terms have also been studied in [4,26]. In this case, we provide sufficient conditions for exponential stability in the spatial sup norm with respect to boundary disturbances (Theorem 2.6). This is an important result for control purposes, because it shows that boundary controllers designed with the backstepping methodology are robust with respect to diffusion (which is a high-order perturbation term). The obtained result is independent of the speed of the hyperbolic PDE and can be interpreted as a delay-independent stability condition for the corresponding parabolic PDE with delayed trace terms.

The present work is structured as follows: the main results of the paper are stated in Section 2. Section 3 contains the proofs of all results. Finally, the concluding remarks of the present work are presented in Section 4.

**Notation.** Throughout this paper, we adopt the following notation.

* $\Re_+ := [0,+\infty)$. Let $U \subseteq \Re^n$ be a set with non-empty interior and let $\Omega \subseteq \Re$ be a set. By $C^0(U)$ (or $C^0(U;\Omega)$), we denote the class of continuous mappings on $U$ (which take values in $\Omega$). By $C^k(U)$ (or $C^k(U;\Omega)$), where $k \geq 1$, we denote the class of continuous functions on $U$, which have continuous derivatives of order $k$ on $U$ (and also take values in $\Omega$).

* $L^2(0,1)$ denotes the equivalence class of measurable functions $f:[0,1] \to \Re$ for which $\|f\|_2 = \left( \int_0^1 |f(x)|^2 dx \right)^{1/2} < +\infty$. $L^2(0,1)$ is a Hilbert space with inner product $\langle f,g \rangle = \int_0^1 f(x)g(x)dx$. $L^\infty(0,1)$ denotes the equivalence class of measurable functions $f:[0,1] \to \Re$ for which $\|f\|_\infty = ess \sup \{ f(x) : x \in (0,1) \} < +\infty$ (a Banach space).

* Let $u: \Re_+ \times [0,1] \to \Re$ be given. We use the notation $u[t]$ to denote the profile at certain $t \geq 0$, i.e., $(u[t])(x) = u(t,x)$ for all $x \in [0,1]$. When $u(t,x)$ is differentiable with respect to $x \in [0,1]$, we use the notation $u'(t,x)$ for the derivative of $u$ with respect to $x \in [0,1]$, i.e., $u'(t,x) = \frac{\partial u}{\partial x}(t,x)$.

* For an integer $k \geq 1$, $H^k(0,1)$ denotes the Sobolev space of functions in $L^2(0,1)$ with all its weak derivatives up to order $k \geq 1$ in $L^2(0,1)$.

## 2. Main Results

*2.A. Movement of Chemicals Underground*

Certain chemicals are released at position $\xi = 0$ and enter the groundwater system. Let $\phi \in (0,1)$ be the porosity of the soil, $v \geq 0$ be the velocity of the bulk movement of the groundwater, $c(t,\xi)$ and $n(t,\xi)$ be the concentration of the chemicals dissolved in water and the sorbed concentration of chemicals in the soil, respectively, at position $\xi \in [0,L]$ (horizontal coordinate) and time $t \geq 0$.

The physical law that allows us to obtain a mathematical model for this process is Fick's law: the rate per unit area per unit time that mass of chemicals crosses a plane section through the flow at position $\xi \in [0,L]$ and time $t \geq 0$ is equal to $-D\frac{\partial c}{\partial \xi}(t,\xi) + \phi v c(t,\xi)$, where $D > 0$ is the diffusion coefficient. Taking into account that the rate of sorption of chemicals in the soil is proportional to $\frac{\partial n}{\partial t}(t,\xi)$, the mass balance for the chemical gives the equation:

$$\frac{\partial c}{\partial t}(t,\xi) = \frac{D}{\phi}\frac{\partial^2 c}{\partial \xi^2}(t,\xi) - v\frac{\partial c}{\partial \xi}(t,\xi) - \frac{1}{\phi}\frac{\partial n}{\partial t}(t,\xi), \text{ for } (t,\xi) \in (0,+\infty) \times (0,L) \quad (2.1)$$

We assume that the concentration of chemicals dissolved in underground water at $\xi = 0$ is time-varying and takes values around a nominal value $c_0 > 0$. At $\xi = L$ the ground meets the sea, where the concentration of chemicals is zero. Therefore, we obtain the boundary conditions:

$$c(t,0) - c_0 - \tilde{d}(t) = c(t,L) = 0, \text{ for all } t \geq 0 \quad (2.2)$$

where $\tilde{d}$ is the variation of the concentration of chemicals dissolved in water at the source ($\xi = 0$).

In order to complete the description of the mathematical model of the process, we need an empirical relation that provides quantitative information about the rate of sorption of chemicals in the soil. The rate of sorption of chemicals in the soil at position $\xi \in [0,L]$ and time $t \geq 0$ has to be a non-decreasing function of $c(t,\xi)$ and a non-increasing function of $n(t,\xi)$. The simplest relation that describes such a dependence is

$$\frac{\partial n}{\partial t}(t,\xi) = a c(t,\xi) - b n(t,\xi), \text{ for } (t,\xi) \in (0,+\infty) \times (0,L) \quad (2.3)$$

where $a, b > 0$ are constants.

Substituting (2.3) into (2.1) and defining

$$u_1(t,z) := c_0^{-1} \exp\left(-\frac{vL\phi}{2D}z\right)\left(c\left(\frac{L^2\phi}{D}t, Lz\right) - c_{eq}(Lz)\right), \quad u_2(t,z) := c_0^{-1} \exp\left(-\frac{vL\phi}{2D}z\right)\left(n\left(\frac{L^2\phi}{D}t, Lz\right) - n_{eq}(Lz)\right),$$

to be the scaled deviations from the nominal concentration profiles

$$c_{eq}(\xi) = \frac{b}{a}n_{eq}(\xi) = c_0 \frac{\exp\left(\frac{\phi v}{D}L\right) - \exp\left(\frac{\phi v}{D}\xi\right)}{\exp\left(\frac{\phi v}{D}L\right) - 1} \text{ for } \xi \in [0,L], \text{ we obtain from (2.1), (2.2) and (2.3) the}$$

following mathematical model of the process:

$$\frac{\partial u_1}{\partial t}(t,z) - \frac{\partial^2 u_1}{\partial z^2}(t,z) + K u_1(t,z) - r\tilde{b} u_2(t,z) = \frac{\partial u_2}{\partial t}(t,z) - \tilde{a} u_1(t,z) + \tilde{b} u_2(t,z) = 0$$

$$\text{for } (t,z) \in (0,+\infty) \times (0,1) \quad (2.4)$$

$$u_1(t,0) - d(t) = u_1(t,1) = 0, \text{ for } t \geq 0 \quad (2.5)$$

where $d(t) := c_0^{-1}\tilde{d}\left(\frac{L^2\phi}{D}t\right)$, $\tilde{a} := a\frac{L^2\phi}{D}$, $\tilde{b} := b\frac{L^2\phi}{D}$, $r := \phi^{-1}$, $K := L^2\frac{v^2\phi^2 + 4aD}{4D^2}$. All parameters and variables appearing in model (2.4), (2.5) are dimensionless.

System (2.4), (2.5) is the feedback interconnection of a parabolic PDE with a first-order zero-speed hyperbolic PDE (or an infinitely-parameterized scalar ODE). Its dynamical behavior is very different from that of a parabolic PDE: to see this notice that system (2.4), (2.5) may be transformed to a wave equation (or Klein-Gordon equation) with Kelvin-Voigt damping that may also include viscous damping and stiffness terms

$$\frac{\partial^2 u_1}{\partial t^2}(t,z) + (\tilde{b}+K)\frac{\partial u_1}{\partial t}(t,z) = \frac{\partial^3 u_1}{\partial z^2 \partial t}(t,z) + \tilde{b}\frac{\partial^2 u_1}{\partial z^2}(t,z) + \tilde{b}(r\tilde{a}-K)u_1(t,z), \text{ for } (t,z) \in (0,+\infty)\times(0,1)$$

and, conversely, any wave (or Klein-Gordon) equation with Kelvin-Voigt damping can be transformed to the form (2.4), (2.5).

We next provide existence/uniqueness results for the initial-boundary value (2.4), (2.5) with

$$u_1[0] = u_{1,0}, \quad u_2[0] = u_{2,0} \tag{2.6}$$

where $u_{1,0}, u_{2,0}$ are real functions on $[0,1]$. Our main result is the following theorem.

**Theorem 2.1 (Existence/Uniqueness):** *Consider the initial-boundary value problem (2.4), (2.5), (2.6), where $K, r, \tilde{a}, \tilde{b} \in \Re$ are constants. For every $u_{2,0} \in C^1([0,1])$, $u_{1,0} \in \{w \in H^3(0,1): w(1) = w''(0) = w''(1) = 0\}$ and for every disturbance input $d \in C^2(\Re_+)$ with $d(0) = u_{1,0}(0)$, there exists a unique pair of mappings $u_1 \in C^0(\Re_+ \times [0,1]) \cap C^1((0,+\infty)\times[0,1])$, $u_2 \in C^1(\Re_+ \times [0,1])$ with $u_1[t] \in C^2([0,1])$ for $t > 0$ satisfying (2.4), (2.5), (2.6).*

Theorem 2.1 implies that there exists a set of initial conditions $S \subseteq \{w \in C^2([0,1]): w(1) = 0\} \times C^1([0,1])$ with the following property:

"For every $u_0 = (u_{1,0}, u_{2,0}) \in S$ there exists a non-empty set $\Phi(u_0) \subseteq \{d \in C^0(\Re_+) \cap C^1((0,+\infty)): d(0) = u_{1,0}(0)\}$ such that for every disturbance input $d \in \Phi(u_0)$ the initial-boundary value problem (2.4), (2.5), (2.6), there exists a unique pair of mappings $u_1 \in C^0(\Re_+ \times [0,1]) \cap C^1((0,+\infty)\times[0,1])$, $u_2 \in C^1(\Re_+ \times [0,1])$ with $(u_1[t], u_2[t]) \in S$ for $t > 0$ satisfying (2.4), (2.5), (2.6)."

Moreover, it holds that $\{w \in H^3(0,1): w(1) = w''(0) = w''(1) = 0\} \times C^1([0,1]) \subseteq S$.

We next provide sufficient conditions for ISS in the sup norm for system (2.4), (2.5). The ISS result is given in the following theorem.

**Theorem 2.2 (ISS in the spatial sup norm):** *Consider the hyperbolic-parabolic system (2.4), (2.5), where $K, r, \tilde{a} \in \Re$, $\tilde{b} > 0$ are constants. Suppose that*

$$|r\tilde{a}| < K + \pi^2. \tag{2.7}$$

*Then there exist constants $M, \gamma, \delta > 0$ such that for every $u_{1,0} \in C^0([0,1])$ with $u_{1,0}(1) = 0$, $u_{2,0} \in C^1([0,1])$ and for every disturbance input $d \in C^0(\Re_+)$ with $d(0) = u_{1,0}(0)$, for which there exists a unique pair of mappings $u_1 \in C^0(\Re_+ \times [0,1]) \cap C^1((0,+\infty)\times[0,1])$, $u_2 \in C^1(\Re_+ \times [0,1])$ with $u_1[t] \in C^2([0,1])$ for $t > 0$ satisfying (2.4), (2.5), (2.6), the following inequality holds for all $t \geq 0$:*

$$\|u_1[t]\|_\infty + \|u_2[t]\|_\infty \leq M \exp(-\delta t)(\|u_{1,0}\|_\infty + \|u_{2,0}\|_\infty) + \gamma \max_{0 \leq s \leq t}(|d(s)|) \tag{2.8}$$

Theorem 2.2 is proved by applying the small-gain methodology in conjunction with the ISS estimates in the spatial sup norm for parabolic PDEs. It should be noticed that the stability condition (2.7) is sharp: when $r\tilde{a} \geq K + \pi^2$ there exist solutions of (2.4), (2.5) with $d(t) \equiv 0$ of the form $u_1(t,z) = \exp(\mu t)\sin(\pi z)$, $u_2(t,z) = k\exp(\mu t)\sin(\pi z)$ for which $\mu \geq 0$, and consequently these solutions do not tend to zero. An interpretation of (2.7) can be made by looking at the wave equation with Kelvin-Voigt damping that corresponds to the hyperbolic-parabolic system (2.4), (2.5), namely, the equation $\frac{\partial^2 u}{\partial t^2}(t,z) + (\tilde{b}+K)\frac{\partial u}{\partial t}(t,z) = \frac{\partial^3 u}{\partial z^2 \partial t}(t,z) + \tilde{b}\frac{\partial^2 u}{\partial z^2}(t,z) + \tilde{b}(r\tilde{a}-K)u(t,z)$. In this case, condition (2.7) implies that a possible anti-stiffness does not dominate the strain. For the case where the stiffness term is absent, we obtain the following corollary.

**Corollary 2.3 (ISS in the spatial sup norm):** *Consider the wave equation with Kelvin-Voigt and viscous damping*

$$\frac{\partial}{\partial t}\left(\frac{\partial u}{\partial t}(t,z) - \sigma\frac{\partial^2 u}{\partial z^2}(t,z)\right) = c^2 \frac{\partial^2 u}{\partial z^2}(t,z) - \mu\frac{\partial u}{\partial t}(t,z), \text{ for } (t,z) \in (0,+\infty) \times (0,1) \quad (2.9)$$

$$u(t,0) - d(t) = u(t,1) = 0, \text{ for } t \geq 0 \quad (2.10)$$

*where* $\sigma, c > 0$, $\mu \geq 0$ *are constants. Suppose that*

$$2c^2 < 2\mu\sigma + \sigma^2 \pi^2. \quad (2.11)$$

*Then there exist constants* $M, \gamma, \delta > 0$ *such that for every* $u_0 \in C^2([0,1])$, $w_0 \in C^0([0,1])$ *with* $u_0(1) = 0$, $(w_0 - \sigma u_0'') \in C^1([0,1])$ *and for every disturbance input* $d \in C^0(\Re_+)$ *with* $d(0) = u_0(0)$, *for which there exists a unique pair of mappings* $u \in C^1(\Re_+ \times [0,1])$, $w \in C^0(\Re_+ \times [0,1])$ *with* $(w - \sigma u'') \in C^1(\Re_+ \times [0,1])$, $u[t] \in C^2([0,1])$ *for* $t > 0$ *satisfying (2.9), (2.10),* $u[0] = u_0$ *and* $\frac{\partial u}{\partial t}[0] = w_0$, *the following inequality holds for all* $t \geq 0$:

$$\|u[t]\|_\infty + \left\|\frac{\partial u}{\partial t}[t] - \sigma u''[t]\right\|_\infty \leq M \exp(-\delta t)(\|u_0\|_\infty + \|w_0 - \sigma u_0''\|_\infty) + \gamma \max_{0 \leq s \leq t}(|d(s)|) \quad (2.12)$$

It is interesting to notice that the coefficient $\gamma > 0$ appearing in the estimate

$$\|u[t]\|_\infty \leq M \exp(-\delta t)(\|u_0\|_\infty + \|w_0 - \sigma u_0''\|_\infty) + \gamma \max_{0 \leq s \leq t}(|d(s)|) \quad (2.13)$$

can be interpreted as the magnification factor of a boundary oscillation to the main body of a string, which has the other end pinned down. Due to the fact that $u(t,0) = d(t)$, the coefficient $\gamma > 0$ appearing in (2.13) is always greater or equal to 1. An estimate of the magnification factor $\gamma > 0$ can be obtained by following the proof of Theorem 2.2. Indeed, it is shown that under assumption (2.11), for every $\theta \in \left(0, \frac{\pi}{2}\right)$ and for every $\varepsilon > 0$ with $\frac{(1+\varepsilon)^2|\sigma\mu - c^2|}{\sigma\mu - c^2 + \sigma^2(\pi - 2\theta)^2} < 1$, there exist constants $M, \delta > 0$ such that the following estimate holds for all $t \geq 0$:

$$\|u[t]\|_\infty \leq M \exp(-\delta t)(\|u_0\|_\infty + \|w_0 - \sigma u_0''\|_\infty) + \frac{1+\varepsilon}{\sin(\theta)(1 - (1+\varepsilon)\sqrt{P(\theta)})^2} \max_{0 \leq s \leq t}(|d(s)|) \quad (2.14)$$

where $P(\theta) := \frac{|\sigma\mu - c^2|}{\sigma\mu - c^2 + \sigma^2(\pi - 2\theta)^2}$. Therefore, it follows that estimate (2.13) holds for all $\gamma > g\left(\frac{c^2 - \mu\sigma}{\sigma^2}\right)$, where $g(s) := \inf\left\{\frac{1}{\sin(\theta)(1-\sqrt{P(\theta)})^2} : 0 < \theta < \frac{\pi - \sqrt{|s|-s}}{2}\right\}$, provided that (2.11) holds. The graph of the function $g$ is shown in Fig.1. It is shown that $g$ has a unique minimum at 0, which implies that the smallest lower bound for the magnification factor $\gamma > 0$ is obtained when $\mu\sigma = c^2$, which is equal to 1.

Turning back to the application that motivated the study of system (2.4), (2.5) and using the definitions $\tilde{a} := a\frac{L^2\phi}{D}$, $\tilde{b} := b\frac{L^2\phi}{D}$, $r := \phi^{-1}$, $K := L^2\frac{v^2\phi^2 + 4aD}{4D^2}$, we conclude that the stability condition (2.7) holds automatically for the problem of the underground movement of chemicals. Therefore, the deviation of the concentration of chemicals (both in water and in soil) satisfy the ISS property with respect to time-varying variations of the concentration of chemicals at the source point.

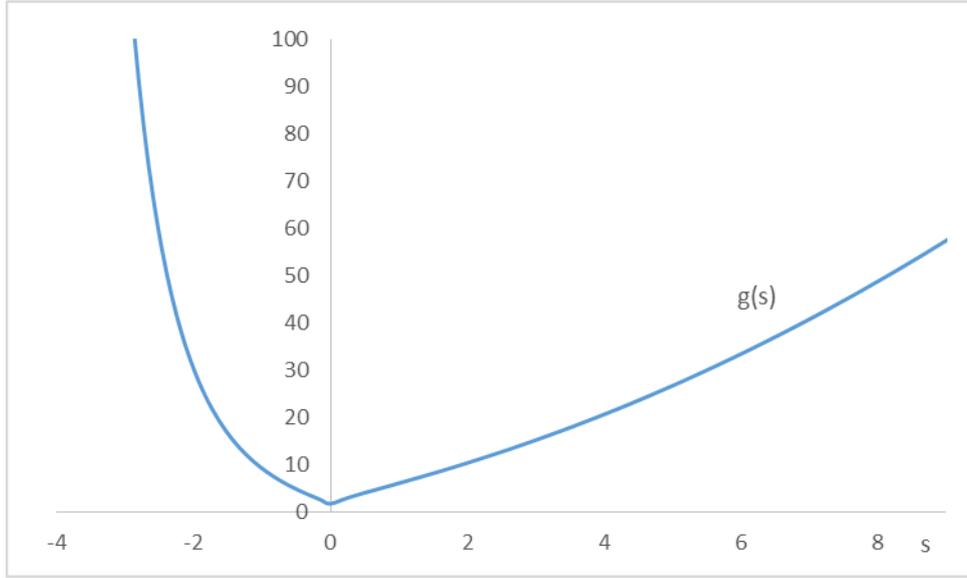

**Figure 1:** The graph of the function $g$.

*2.B. Boundary and In Domain Feedback Interconnection*

In this section we present the analysis of a parabolic PDE which is interconnected with a first-order hyperbolic PDE by means of a combination of boundary and in-domain terms. The overall system contains a non-local term and is analyzed both with respect to existence/uniqueness of solutions as well as with respect to exponential stability in the spatial sup norm.

Consider the following system of PDEs

$$\frac{\partial u_1}{\partial t}(t,z) - p\frac{\partial^2 u_1}{\partial z^2}(t,z) - au_1(t,z) - \int_0^1 b(z,s)u_2(t,s)ds = \frac{\partial u_2}{\partial t}(t,z) + c\frac{\partial u_2}{\partial z}(t,z) = 0$$

for $(t,z) \in (0,+\infty) \times (0,1)$ (2.15)

$$u_1(t,0) = \frac{\partial u_1}{\partial z}(t,1) - qu_1(t,1) = 0, \text{ for } t \geq 0 \quad (2.16)$$

$$u_2(t,0) = ku_1(t,1), \text{ for } t \geq 0 \quad (2.17)$$

where $p, c > 0$, $q < 1$, $a, k \in \Re$ are constants and $b \in C^1([0,1]^2)$ is a given function. Model (2.15), (2.16), (2.17) arises in the study of liquid metal droplet-generated "extreme ultraviolet" light sources for photolithography. It is a system of a parabolic PDE which is interconnected with a first-order hyperbolic PDE, by means of two different terms: the in-domain, non-local term $\int_0^1 b(z,s)u_2(t,s)ds$ that appears in the parabolic PDE and the boundary non-local trace term $ku_1(t,1)$ that appears in the boundary condition (2.17).

Model (2.15), (2.16), (2.17) is also closely related to a parabolic PDE with delays. Indeed, solving the hyperbolic PDE in (2.15), we obtain the following parabolic PDE with distributed delays:

$$\frac{\partial u_1}{\partial t}(t,z) = p\frac{\partial^2 u_1}{\partial z^2}(t,z) + au_1(t,z) + ck\int_{t-c^{-1}}^{t} b(z,c(t-s))u_1(s,1)ds$$

In the above setting, we are interested in obtaining delay-independent stability conditions that guarantee exponential stability for any value of the delay, i.e., stability conditions that are independent of $c > 0$.

Let $A: D \to L^2(0,1)$ be the SL operator associated with the parabolic PDE in (2.15), namely, the linear operator $(Au)(z) = -pu''(z) - au(z)$ for $z \in [0,1]$, defined on

$$D := \left\{ u \in H^2(0,1) : u(0) = u'(1) - qu(1) = 0 \right\} \tag{2.18}$$

The eigenfunctions of $A$ are expressed by $\phi_n(z) = A_n \sin\left(\left((2n-1)\frac{\pi}{2} - b_n\right)z\right)$ for $n = 1,2,\ldots$, $z \in [0,1]$, where $b_n \in \left(-\frac{\pi}{2}, \frac{\pi}{2}\right)$, $n = 1,2,\ldots$ are the unique solutions on $\left(-\frac{\pi}{2}, \frac{\pi}{2}\right)$ of the equations $\left((2n-1)\frac{\pi}{2} - b_n\right)\cot\left((2n-1)\frac{\pi}{2} - b_n\right) = q$ and $A_n = \sqrt{\dfrac{2(2n-1)\pi - 4b_n}{(2n-1)\pi - 2b_n - \sin((2n-1)\pi - 2b_n)}}$. The eigenvalues of $A$ are given by $\lambda_n = p\left((2n-1)\frac{\pi}{2} - b_n\right)^2 - a$, for $n = 1,2,\ldots$. We assume that all eigenvalues of $A$ are positive, i.e.,

$$p(\pi - 2b_1)^2 > 4a \tag{2.19}$$

We next provide existence/uniqueness results for the initial-boundary value (2.15), (2.16), (2.17) with

$$u_1[0] = u_{1,0}, \quad u_2[0] = u_{2,0} \tag{2.20}$$

where $u_{1,0}, u_{2,0}$ are real functions on $[0,1]$. Our main result is the following theorem.

**Theorem 2.4 (Existence/Uniqueness):** *Consider the initial-boundary value problem (2.15), (2.16), (2.17), (2.20), where $p, c > 0$, $q < 1$, $a, k \in \Re$ are constants and $b \in C^1([0,1]^2)$ is a given function and suppose that (2.19) holds. Then for every $u_{1,0} \in D$, $u_{2,0} \in C^1([0,1])$ with $u_{2,0}(0) = k u_{1,0}(1)$, $\theta_{1,0} \in D$ $u'_{2,0}(0) = -c^{-1} k \theta_{1,0}(1)$, where $\theta_{1,0}(z) = p u''_{1,0}(z) + a u_{1,0}(z) + \int_0^1 b(z,s) u_{2,0}(s) ds$ for $z \in [0,1]$, there exists a unique pair of mappings $u_1 \in C^0(\Re_+ \times [0,1]) \cap C^1((0,+\infty) \times [0,1])$, $u_2 \in C^1(\Re_+ \times [0,1])$ with $u_1[t] \in C^2([0,1])$ for $t > 0$ satisfying (2.15), (2.16), (2.17), (2.20).*

We next provide sufficient conditions for exponential stability in the sup norm for system (2.15), (2.16), (2.17), which are independent of $c > 0$. We first need the following auxiliary lemma.

**Lemma 2.5:** *Suppose that (2.19) holds. Then there exist constants $\theta \in (0, \pi)$, $\omega \in [0, \pi - \theta)$ with $\omega \cot(\omega + \theta) > q$ and $p\omega^2 > a$.*

**Theorem 2.6 (Exponential Stability in the sup norm independent of $c > 0$):** *Consider the hyperbolic-parabolic system (2.15), (2.16), (2.17), where $p, c > 0$, $q < 1$, $a, k \in \Re$ are constants and $b \in C^1([0,1]^2)$ is a given function and suppose that (2.19) holds. Moreover, suppose that*

$$|k| \max_{0 \le z \le 1}\left(\frac{\sin(\omega + \theta)}{\sin(\omega z + \theta)} \int_0^1 |b(z,s)| ds\right) < p\omega^2 - a \tag{2.21}$$

*for certain constants $\theta \in (0, \pi)$, $\omega \in [0, \pi - \theta)$ with $p\omega^2 > a$ and $\omega \cot(\omega + \theta) > q$, whose existence is established by Lemma 2.5. Then there exist constants $M, \delta > 0$ such that for every $u_{1,0}, u_{2,0} \in C^1([0,1])$ with $u_{1,0}(1) = u'_{1,0}(1) - q u_{1,0}(1) = u_{2,0}(0) - k u_{1,0}(1) = 0$, for which there exists a unique pair of mappings $u_1 \in C^0(\Re_+ \times [0,1]) \cap C^1((0,+\infty) \times [0,1])$, $u_2 \in C^1(\Re_+ \times [0,1])$ with $u_1[t] \in C^2([0,1])$ for $t > 0$ satisfying (2.15), (2.16), (2.17), (2.20), the following inequality holds for all $t \ge 0$:*

$$\|u_1[t]\|_\infty + \|u_2[t]\|_\infty \le M \exp(-\delta t)\left(\|u_{1,0}\|_\infty + \|u_{2,0}\|_\infty\right) \tag{2.22}$$

Inequality (2.21) imposes a bound on the product of the static gains of the interconnecting, nonlocal terms $\int_0^1 b(z,s)u_2(t,s)ds$, $ku_1(t,1)$ that may lead the solution far from equilibrium. Inequality (2.21) is independent of $c>0$ and consequently, may be a conservative stability condition. In the context of the related delay parabolic PDE $\frac{\partial u_1}{\partial t}(t,z) = p\frac{\partial^2 u_1}{\partial z^2}(t,z) + a u_1(t,z) + ck\int_{t-c^{-1}}^{t} b(z,c(t-s))u_1(s,1)ds$, condition (2.21) is a delay-independent stability condition. Delay-independent stability conditions are conservative conditions which are useful in practice because the delays are not easily measured or may vary.

Another use of the stability condition (2.21) can be shown by studying a crucial question that often arises in engineering practice: is it safe to ignore diffusion? The following example shows how the engineer can exploit the stability condition (2.21) in order to answer the question regarding robustness to diffusion.

**Example (Is It Safe to Ignore Diffusion?):** Applying backstepping to a pair of hyperbolic PDEs (see [6,10]), we obtain the following hyperbolic PDE-PDE loop

$$\frac{\partial w_1}{\partial t}(t,z) + v\frac{\partial w_1}{\partial z}(t,z) = \frac{\partial w_2}{\partial t}(t,z) + c\frac{\partial w_2}{\partial z}(t,z) = 0, \text{ for } (t,z) \in (0,+\infty) \times (0,1) \quad (2.23)$$

$$w_1(t,0) = w_2(t,0) - kw_1(t,1) = 0, \text{ for } t \geq 0 \quad (2.24)$$

where $v,c > 0$, $k \in \Re$ are constants. The equilibrium point $0 \in \left(C^0([0,1])\right)^2$ is exponentially stable (in fact, finite-time stable). However, when diffusion phenomena are present in one of the PDEs, then the actual closed-loop system is described by (2.24) along with the equations

$$\frac{\partial w_1}{\partial t}(t,z) - p\frac{\partial^2 w_1}{\partial z^2}(t,z) + v\frac{\partial w_1}{\partial z}(t,z) - p\int_0^1 l(z,s)w_2(t,s)ds = \frac{\partial w_2}{\partial t}(t,z) + c\frac{\partial w_2}{\partial z}(t,z) = 0$$

$$\text{for } (t,z) \in (0,+\infty) \times (0,1) \quad (2.25)$$

$$\frac{\partial w_1}{\partial z}(t,1) = 0, \text{ for } t \geq 0 \quad (2.26)$$

where $p > 0$ is a constant and $l \in C^1\left([0,1]^2\right)$ is a given function. Does exponential stability for system (2.23), (2.24) guarantee exponential stability for system (2.24), (2.25), (2.26) when the diffusion coefficient $p > 0$ is sufficiently small?

The answer is "yes". To prove this, we perform the transformation

$$u_1(t,z) = \exp\left(-\frac{v}{2p}z\right)w_1(t,z) \quad , \quad u_2(t,z) = w_2(t,z), \text{ for } t \geq 0, \; z \in [0,1] \quad (2.27)$$

which transforms system (2.24), (2.25), (2.26) to system (2.15), (2.16), (2.17), with

$$q := -\frac{v}{2p}, \; b(z,s) = p\exp\left(-\frac{v}{2p}z\right)l(z,s) \text{ for } z,s \in [0,1], \; a := -\frac{v^2}{4p}$$

Applying Theorem 2.6 with $\theta = \pi/2$, $\omega = 0$, we conclude that $0 \in \left(C^0([0,1])\right)^2$ is exponentially stable in the sup norm for system (2.24), (2.25), (2.26), provided that the following condition holds

$$2p\sqrt{|k|\max_{0\leq z\leq 1}\left(\int_0^1 |l(z,s)|ds\right)} < v \quad (2.28)$$

Consequently, it follows from (2.28) that there exists $P \in (0,+\infty]$ such that $0 \in \left(C^0([0,1])\right)^2$ is exponentially stable in the sup norm for system (2.24), (2.25), (2.26), provided that $p \in (0,P)$. Notice that inequality (2.28) provides an explicit lower bound for $P \in (0,+\infty]$, which is independent of $c > 0$. ◁

# 3. Proofs of Main Results

The proofs of the main results require two auxiliary technical results that can be obtained by modifying slightly the proofs of existing results in the literature.

The first auxiliary technical result deals with the Sturm-Liouville (SL) operator $A: D \to L^2(0,1)$ defined by

$$(Af)(z) = -\frac{d}{dz}\left(p(z)\frac{df}{dz}(z)\right) + q(z)f(z), \text{ for all } f \in D \text{ and } z \in (0,1) \quad (3.1)$$

where $p \in C^1([0,1];(0,+\infty))$, $q \in C^0([0,1];\Re)$ and $D \subseteq H^2(0,1)$ is given by

$$D := \{f \in H^2(0,1): b_1 f(0) + b_2 f'(0) = a_1 f(1) + a_2 f'(1) = 0\} \quad (3.2)$$

where $a_1, a_2, b_1, b_2$ are real constants with $|a_1| + |a_2| > 0$, $|b_1| + |b_2| > 0$. It is well-known that all eigenvalues of the SL operator $A: D \to L^2(0,1)$, defined by (3.1), (3.2) are real. The eigenvalues form an infinite, increasing sequence $\lambda_1 < \lambda_2 < \ldots < \lambda_n < \ldots$ with $\lim_{n \to \infty}(\lambda_n) = +\infty$ and to each eigenvalue $\lambda_n \in \Re$ ($n = 1, 2, \ldots$) corresponds exactly one eigenfunction $\phi_n \in C^2([0,1]; \Re)$ that satisfies $A\phi_n = \lambda_n \phi_n$ and $b_1 \phi_n(0) + b_2 \phi_n'(0) = a_1 \phi_n(1) + a_2 \phi_n'(1) = 0$. Moreover, the eigenfunctions form an orthonormal basis of $L^2(0,1)$. We use the following assumptions for the SL operator $A: D \to L^2(0,1)$ defined by (3.1), (3.2), where $a_1, a_2, b_1, b_2$ are real constants with $|a_1| + |a_2| > 0$, $|b_1| + |b_2| > 0$.

**(H1)** *Either $b_2 > 0$ or $b_2 = 0$ and $b_1 < 0$,*
**(H2)** *Either $a_2 > 0$ or $a_2 = 0$ and $a_1 > 0$.*
**(H3):** *The eigenvalues and eigenfunctions of the SL operator $A: D \to L^2(0,1)$ defined by (3.1), (3.2), satisfy*

$$\sum_{n=N}^{\infty} \lambda_n^{-1} \max_{0 \leq z \leq 1}(|\phi_n(z)|) < +\infty, \text{ for certain } N > 0 \text{ with } \lambda_N > 0 \quad (3.3)$$

It is important to notice that the validity of Assumption (H3) can be verified without the knowledge of the eigenvalues and eigenfunctions of the SL operator $A$. More specifically, it is shown in [25] that assumption (H3) holds automatically, if $b_2, a_1, a_2 \geq 0$, $b_1 \leq 0$ and the function $p$ is of class $C^2([0,1];(0,+\infty))$.

A careful inspection of the proof of Theorem 3.1 in [13] allows us to obtain the following useful corollary, which is the first auxiliary technical result needed in the proofs of the main results.

**Corollary 4.1:** *Consider the SL operator $A: D \to L^2(0,1)$ defined by (3.1), (3.2), under Assumptions (H1), (H2), (H3). Let $T > 0$ be a constant and let $f \in C^0([0,T] \times [0,1])$ be a given function for which $(0,T] \times [0,1] \ni (t,z) \to \frac{\partial f}{\partial t}(t,z)$ is continuous, with $f[t]$ being a piecewise $C^1$ function on $[0,1]$ for all $t \in [0,T]$. Then for every $u_0 \in D$ the function $u: [0,T] \times [0,1] \to \Re$ defined for $(t,z) \in [0,T] \times [0,1]$ by*

$$u(t,z) = \sum_{n=1}^{\infty} \phi_n(z) \exp(-\lambda_n t) \int_0^1 \phi_n(s) u_0(s) ds + \sum_{n=1}^{\infty} \phi_n(z) \int_0^t \exp(-\lambda_n(t-\tau)) \left(\int_0^1 \phi_n(s) f(\tau,s) ds\right) d\tau \quad (3.4)$$

*is of class $u \in C^0([0,T] \times [0,1]) \cap C^1((0,T] \times [0,1])$ satisfying $u[t] \in C^2([0,1])$ for all $t \in (0,T]$, $u[0] = u_0$ and*

$$\frac{\partial u}{\partial t}(t,z) - \frac{\partial}{\partial z}\left(p(z)\frac{\partial u}{\partial z}(t,z)\right) + q(z)u(t,z) = f(t,z), \text{ for all } (t,z) \in (0,T] \times (0,1) \quad (3.5)$$

$$b_1 u(t,0) + b_2 \frac{\partial u}{\partial z}(t,0) = a_1 u(t,1) + a_2 \frac{\partial u}{\partial z}(t,1) = 0, \text{ for all } t \in (0,T] \quad (3.6)$$

The second auxiliary result deals with the parabolic PDE (3.5) with boundary conditions given by

$$b_1 u(t,0) + b_2 \frac{\partial u}{\partial z}(t,0) - d_0(t) = a_1 u(t,1) + a_2 \frac{\partial u}{\partial z}(t,1) - d_1(t) = 0, \quad (3.7)$$

where $f(t,z), d_0(t), d_1(t)$ are external inputs. We work with inputs that belong to a set that possibly depends on the initial condition, as specified by the following definition.

**Definition 4.2:** *Consider the SL operator $A: D \to L^2(0,1)$ defined by (3.1), (3.2), under Assumptions (H1), (H2), (H3). For every given $u_0 \in H^2(0,1)$, $\tilde{\Phi}(A; u_0)$ denotes the non-empty set of disturbance inputs for which the following implication holds:*
*"If $(f, d_0, d_1) \in \tilde{\Phi}(A; u_0)$ then $d_0, d_1 \in C^0(\Re_+)$, $f \in C^0(\Re_+ \times [0,1])$ and the evolution equation (3.5) with (3.7) and initial condition $u_0 \in H^2(0,1)$ has a unique solution $u \in C^0(\Re_+ \times [0,1])$ with $u \in C^1((0,+\infty) \times [0,1])$ satisfying $u[t] \in C^2([0,1])$ for all $t > 0$, $u[0] = u_0$, (3.5) for all $(t, z) \in (0,+\infty) \times (0,1)$ and (3.7) for all $t \geq 0$".*

Theorem 2.1 in [14] and Corollary 4.1 guarantee that if $f \in C^0(\Re_+ \times [0,1])$ is a function for which $(0,+\infty) \times [0,1] \ni (t,z) \to \frac{\partial f}{\partial t}(t,z)$ is continuous, with $f[t]$ being a piecewise $C^1$ function on $[0,1]$ for all $t \geq 0$, and if $d_0, d_1 \in C^2(\Re_+)$ satisfy $b_1 u_0(0) + b_2 u_0'(0) - d_0(0) = a_1 u_0(1) + a_2 u_0'(1) - d_1(0) = 0$ then $(f, d_0, d_1) \in \tilde{\Phi}(A; u_0)$. Notice that if $(f, d_0, d_1) \in \tilde{\Phi}(A; u)$ and $(\tilde{f}, \tilde{d}_0, \tilde{d}_1) \in \tilde{\Phi}(A; v)$ for certain $u, v \in D$ then $(\lambda f, \lambda d_0, \lambda d_1) \in \tilde{\Phi}(A; \lambda u)$ for every $\lambda \in \Re$ and $(f + \tilde{f}, d_0 + \tilde{d}_0, d_1 + \tilde{d}_1) \in \tilde{\Phi}(A; u+v)$.

In order to derive ISS estimates expressed in the spatial $L^\infty$ norm of the solution of (3.5), (3.7), we need the following assumption.

**(H4)** *There exists a function $\eta \in C^2([0,1]; (0,+\infty))$ and a constant $\sigma > 0$ such that $p(z)\eta''(z) + p'(z)\eta'(z) - q(z)\eta(z) \leq -\sigma r(z)\eta(z)$ for all $z \in [0,1]$. Moreover, the inequalities $b_1 \eta(0) + b_2 \eta'(0) < 0$ and $a_1 \eta(1) + a_2 \eta'(1) > 0$ hold.*

Assumption (H4) is strongly related to the existence of a positive eigenfunction of the Sturm-Liouville operator $A$ defined by (3.1), (3.2). The problem of the existence of positive eigenfunctions for elliptic operators has been studied in the literature (see [27] on page 112 and references therein). However, it must be noted that here there is a degree of freedom: no specific boundary conditions are assumed to hold for the function $\eta \in C^2([0,1]; (0,+\infty))$: the function is only required to satisfy the inequalities $b_1 \eta(0) + b_2 \eta'(0) < 0$ and $a_1 \eta(1) + a_2 \eta'(1) > 0$.

A careful inspection of the proof of Theorem 2.2 in [14] allows us to obtain the following corollary, which is the second auxiliary technical result needed for the proof of the main results.

**Corollary 4.3 (ISS in the sup-norm):** *Consider the SL operator $A: D \to L^2(0,1)$ defined by (3.1), (3.2), under Assumptions (H1), (H2), (H3), (H4). Then for every $u_0 \in H^2(0,1)$, $(f, d_0, d_1) \in \tilde{\Phi}(A; u_0)$, the unique solution $u \in C^0(\Re_+ \times [0,1]) \cap C^1((0,+\infty) \times [0,1])$ of the evolution equation (3.5) with (3.7) and initial condition $u_0 \in H^2(0,1)$ satisfies the following estimate for all $t \geq 0$:*

$$\|u[t]\|_{\infty,\eta} \leq \max\left(\exp(-\sigma t)\|u_0\|_{\infty,\eta}, \frac{\max_{0 \leq s \leq t}(|d_0(s)|)}{|b_1 \eta(0) + b_2 \eta'(0)|}, \frac{\max_{0 \leq s \leq t}(|d_1(s)|)}{a_1 \eta(1) + a_2 \eta'(1)}\right) + \sigma^{-1} \max_{0 \leq s \leq t}(\|f[s]\|_{\infty,\eta}) \quad (3.8)$$

*where*

$$\|u[t]\|_{\infty,\eta} := \max_{0 \leq z \leq 1}\left(\frac{|u(t,z)|}{\eta(z)}\right). \quad (3.9)$$

We are now ready to give the proofs of the main results of the present work.

**Proof of Theorem 2.1:** It suffices to show that for every $T > 0$, $u_{1,0} \in \{w \in H^3(0,1) : w(1) = w''(0) = w''(1) = 0\}$, $u_{2,0} \in C^1([0,1])$ and $d \in C^2(\Re_+)$ with $d(0) = u_{1,0}(0)$ there exists a unique pair of mappings $u_1 \in C^0([0,T] \times [0,1]) \cap C^1((0,T] \times [0,1])$, $u_2 \in C^1([0,T] \times [0,1])$ with $u_1[t] \in C^2([0,1])$ for $t \in (0,T]$ satisfying (2.4), (2.5), (2.6). Let $T > 0$, $u_{1,0} \in \{w \in H^3(0,1) : w(1) = w''(0) = w''(1) = 0\}$, $u_{2,0} \in C^1([0,1])$ and $d \in C^2(\Re_+)$ with $d(0) = u_{1,0}(0)$ be given (arbitrary).

We first perform the homogenization of the boundary conditions. We define for $(t,z) \in \Re_+ \times [0,1]$:

$$w_1(t,z) := u_1(t,z) - (1-z)d(t) \quad , \quad w_2(t,z) := u_2(t,z)$$
$$f_1(t,z) := -(1-z)\left(\dot{d}(t) + Kd(t)\right) \quad , \quad f_2(t,z) := \tilde{a}(1-z)d(t) \tag{3.10}$$

Using (2.4), (2.5), (2.6) and (3.10), we obtain the equivalent initial-boundary value problem:

$$\frac{\partial w_1}{\partial t}(t,z) - \frac{\partial^2 w_1}{\partial z^2}(t,z) + Kw_1(t,z) - r\tilde{b} w_2(t,z) - f_1(t,z) = \frac{\partial w_2}{\partial t}(t,z) - \tilde{a}w_1(t,z) + \tilde{b}w_2(t,z) - f_2(t,z) = 0$$

$$\text{for } (t,z) \in (0,T] \times (0,1) \tag{3.11}$$

$$w_1(t,0) = w_1(t,1) = 0, \text{ for } t \in [0,T] \tag{3.12}$$

$$w_1[0] = w_{1,0}, \quad w_2[0] = w_{2,0} \tag{3.13}$$

where $w_{1,0}(z) = u_{1,0}(z) - (1-z)d(0)$ for $z \in [0,1]$ and $w_{2,0} = u_{2,0}$. It follows from Corollary 4.1 (with direct computation of the eigenvalues and eigenfunctions of the SL operator that corresponds to the parabolic PDE) and integration of the hyperbolic PDE that any classical solution of (3.11), (3.12), (3.13), $w_1 \in C^0([0,T] \times [0,1]) \cap C^1((0,T] \times [0,1])$, $w_2 \in C^1([0,T] \times [0,1])$ with $w_1[t] \in C^2([0,1])$ for $t \in (0,T]$ and $w_{2,0} \in C^1([0,1])$, $w_{1,0} \in \{w \in H^2(0,1) : w(0) = w(1) = 0\}$, must satisfy the following integral equations for all $(t,z) \in [0,T] \times [0,1]$:

$$w_1(t,z) = 2\sum_{n=1}^{\infty} \sin(n\pi z) \exp\left(-(K + n^2\pi^2)t\right) \int_0^1 \sin(n\pi s) w_{1,0}(s) ds$$

$$+ 2r\tilde{b} \sum_{n=1}^{\infty} \sin(n\pi z) \int_0^t \exp\left(-(K + n^2\pi^2)(t - \tau)\right) \left(\int_0^1 \sin(n\pi s) w_2(\tau, s) ds\right) d\tau \tag{3.14}$$

$$+ 2\sum_{n=1}^{\infty} \sin(n\pi z) \int_0^t \exp\left(-(K + n^2\pi^2)(t - \tau)\right) \left(\int_0^1 \sin(n\pi s) f_1(\tau, s) ds\right) d\tau$$

$$w_2(t,z) = \exp\left(-\tilde{b} t\right) w_{2,0}(z) + \int_0^t \exp\left(-\tilde{b}(t - \tau)\right)\left(\tilde{a} w_1(\tau, z) + f_2(\tau, z)\right) d\tau \tag{3.15}$$

For sufficiently large $k > 0$, the mapping $\left(C^0([0,T] \times [0,1])\right)^2 \ni (v_1, v_2) \to (R_1(v_1, v_2), R_2(v_1, v_2)) \in \left(C^0([0,T] \times [0,1])\right)^2$ defined by

$$(R_1(v_1, v_2))(t,z) = 2\sum_{n=1}^{\infty} \sin(n\pi z) \exp\left(-(K + n^2\pi^2 + k)t\right) \int_0^1 \sin(n\pi s) w_{1,0}(s) ds$$

$$+ 2r\tilde{b} \sum_{n=1}^{\infty} \sin(n\pi z) \int_0^t \exp\left(-(K + n^2\pi^2 + k)(t - \tau)\right) \left(\int_0^1 \sin(n\pi s) v_2(\tau, s) ds\right) d\tau \tag{3.16}$$

$$+ 2\sum_{n=1}^{\infty} \sin(n\pi z) \int_0^t \exp\left(-(K + n^2\pi^2 + k)(t - \tau)\right) \left(\int_0^1 \sin(n\pi s) \exp(-k\tau) f_1(\tau, s) ds\right) d\tau$$

$$(R_2(v_1, v_2))(t,z) = \exp\left(-(k + \tilde{b})t\right) w_{2,0}(z) + \int_0^t \exp\left(-(\tilde{b} + k)(t - \tau)\right)\left(\tilde{a} v_1(\tau, z) + \exp(-k\tau) f_2(\tau, z)\right) d\tau \tag{3.17}$$

is a contraction. Therefore, there exists a unique $(v_1, v_2) \in \left(C^0([0,T]\times[0,1])\right)^2$ satisfying the equation $(v_1, v_2) = (R_1(v_1, v_2), R_2(v_1, v_2))$. By defining $w_i(t,z) = \exp(kt) v_i(t,z)$ for $(t,z) \in [0,T]\times[0,1]$, it follows from (3.16), (3.17) that there exists a unique $(w_1, w_2) \in \left(C^0([0,T]\times[0,1])\right)^2$ satisfying (3.14), (3.15) for $i=1,2$, $(t,z)\in[0,T]\times[0,1]$. Thus, problem (3.11), (3.12), (3.13) has at most one classical solution.

Next, we construct the unique classical solution of the initial-boundary value problem (3.11), (3.12), (3.13). Consider the integral equations for $(t,z)\in[0,T]\times[0,1]$:

$$p_1(t,z) = g(t,z) + 2r\tilde{b} \sum_{n=1}^{\infty} \cos(n\pi z) \int_0^t \exp\left(-(K+n^2\pi^2)(t-\tau)\right) \left(\int_0^1 \left(\cos(n\pi s) - (-1)^n\right) p_2(\tau,s) ds\right) d\tau \quad (3.18)$$

$$p_2(t,z) = \exp\left(-\tilde{b}\, t\right) w'_{2,0}(z) + \tilde{a} \int_0^t \exp\left(-\tilde{b}\,(t-\tau)\right)\left(p_1(\tau,z) - d(\tau)\right) d\tau \quad (3.19)$$

where

$$\begin{aligned}
g(t,z) = &\; 2\sum_{n=1}^{\infty} \cos(n\pi z) \exp\left(-(K+n^2\pi^2)t\right) \int_0^1 \cos(n\pi s) w'_{1,0}(s) ds \\
&- 2r\tilde{b}\, w_{2,0}(0) \sum_{n=1}^{\infty} \left((-1)^n - 1\right)\cos(n\pi z) \int_0^t \exp\left(-(K+n^2\pi^2)(t-\tau)\right) \exp\left(-\tilde{b}\,\tau\right) d\tau \\
&- 2r\tilde{b}\,\tilde{a} \sum_{n=1}^{\infty} \left((-1)^n - 1\right)\cos(n\pi z)\int_0^t \exp\left(-(K+n^2\pi^2)(t-\tau)\right)\left(\int_0^\tau \exp\left(-\tilde{b}\,(\tau-l)\right) d(l)\, dl\right) d\tau \\
&- 2\sum_{n=1}^{\infty} \cos(n\pi z) \int_0^t \exp\left(-(K+n^2\pi^2)(t-\tau)\right)\left(\dot{d}(\tau) + K d(\tau)\right) d\tau
\end{aligned} \quad (3.20)$$

Since $w_{2,0}\in C^1([0,1])$, $w_{1,0}\in \{w\in H^3(0,1): w(0) = w(1) = w''(0) = w''(1) = 0\}$, it follows that for sufficiently large $k>0$, the mapping $\left(C^0([0,T]\times[0,1])\right)^2 \ni (q_1,q_2) \to (\tilde{R}_1(q_1,q_2), \tilde{R}_2(q_1,q_2)) \in \left(C^0([0,T]\times[0,1])\right)^2$ defined by

$$\begin{aligned}
(\tilde{R}_1(q_1,q_2))(t,z) = &\; \exp(-kt) g(t,z) \\
&+ 2r\tilde{b} \sum_{n=1}^{\infty} \cos(n\pi z) \int_0^t \exp\left(-(K+n^2\pi^2+k)(t-\tau)\right)\left(\int_0^1 \left(\cos(n\pi s) - (-1)^n\right) q_2(\tau,s) ds\right) d\tau
\end{aligned} \quad (3.21)$$

$$\begin{aligned}
(\tilde{R}_2(q_1,q_2))(t,z) = &\; -\tilde{a}\exp(-kt) \int_0^t \exp\left(-(k+\tilde{b})(t-\tau)\right) d(\tau) d\tau \\
&+ \exp\left(-(k+\tilde{b})t\right) w'_{2,0}(z) + \tilde{a}\int_0^t \exp\left(-(k+\tilde{b})(t-\tau)\right) q_1(\tau,z) d\tau
\end{aligned} \quad (3.22)$$

is a contraction. Therefore, there exists a unique $(q_1, q_2) \in \left(C^0([0,T]\times[0,1])\right)^2$ satisfying the equation $(q_1, q_2) = (\tilde{R}_1(q_1,q_2), \tilde{R}_2(q_1,q_2))$. By defining $p_i(t,z) = \exp(kt) q_i(t,z)$ for $i=1,2$, $(t,z)\in[0,T]\times[0,1]$, it follows from (3.21), (3.22) that there exists a unique $(p_1, p_2) \in \left(C^0([0,T]\times[0,1])\right)^2$ satisfying (3.18), (3.19) for $(t,z)\in[0,T]\times[0,1]$. It is a matter of straightforward calculations to verify (by using (3.10), (3.18), (3.19)) that the functions $w_1(t,z) = \int_0^z p_1(t,s) ds$ and

$w_2(t,z) = \exp\left(-\tilde{b} t\right) w_{2,0}(0) + \tilde{a}\int_0^t \exp\left(-\tilde{b}(t-\tau)\right) d(\tau) d\tau + \int_0^z p_2(t,s) ds$ satisfy (3.14), (3.15) for $(t,z)\in[0,T]\times[0,1]$.

Notice that $\dfrac{\partial w_1}{\partial z}(t,z) = p_1(t,z)$ is continuous on $[0,T]\times[0,1]$. Consequently, (3.15) implies that $w_2 \in C^1([0,T] \times [0,1])$. Using Corollary 4.1 and (3.14) we obtain that $w_1 \in C^0([0,T] \times [0,1]) \cap C^1((0,T] \times [0,1])$ with $w_1[t] \in C^2([0,1])$ for $t \in (0,T]$. Moreover, (3.11), (3.12), (3.13) hold. The proof is complete. ◁

We next provide the proof of Theorem 2.2.

**Proof of Theorem 2.2:** By virtue of (2.7) there exist $\varepsilon, \zeta > 0$, $\theta \in (0, \pi/2)$ sufficiently small so that

$$L := \frac{(1+\zeta)|r\tilde{a}|(1+\varepsilon)^2}{K + (\pi - 2\theta)^2} < 1 \tag{3.23}$$

Let $u_{1,0} \in C^0([0,1])$ with $u_{1,0}(1) = 0$, $u_{2,0} \in C^1([0,1])$ be (arbitrary) initial conditions and let $d \in C^0(\mathfrak{R}_+)$ with $d(0) = u_{1,0}(0)$ be a disturbance input (arbitrary), for which there exists a unique pair of mappings $u_1 \in C^0(\mathfrak{R}_+ \times [0,1]) \cap C^1((0,+\infty) \times [0,1])$, $u_2 \in C^1(\mathfrak{R}_+ \times [0,1])$ with $u_1[t] \in C^2([0,1])$ for $t > 0$ satisfying (2.4), (2.5), (2.6).

Define the positive function $\eta(z) := \sin(\theta + (\pi - 2\theta)z)$ for $z \in [0,1]$ and the norm $\|u\|_{\infty,\eta} := \max_{0 \le z \le 1}\left(\frac{|u(z)|}{\eta(z)}\right)$. Notice that Assumptions (H1), (H2), (H3), (H4) hold for the PDE problem $\frac{\partial u_1}{\partial t}(t,z) = \frac{\partial^2 u_1}{\partial z^2}(t,z) - K u_1(t,z) + f(t,z)$ with $u_1(t,0) - d(t) = u_1(t,1) = 0$. More specifically, Assumption (H4) holds with $\eta$ as defined above and $\sigma = K + (\pi - 2\theta)^2$.

Applying Corollary 4.3, using the fact that $a + b \le \max((1+\zeta^{-1})a, (1+\zeta)b)$ for all $a,b \ge 0$ and recognizing that the solution of (2.4), (2.5) satisfies the equation $\frac{\partial u_1}{\partial t}(t,z) = \frac{\partial^2 u_1}{\partial z^2}(t,z) - K u_1(t,z) + f(t,z)$ with $f(t,z) := r\tilde{b} u_2(t,z)$, we obtain the following estimate for $t \ge 0$:

$$\|u_1[t]\|_{\infty,\eta} \le \exp(-(K + (\pi - 2\theta)^2)t)\|u_{1,0}\|_{\infty,\eta} + \max_{0 \le s \le t}\left(\frac{(1+\zeta^{-1})|d(s)|}{\sin(\theta)} + \frac{|r|\tilde{b}(1+\zeta)}{K + (\pi - 2\theta)^2}\|u_2[s]\|_{\infty,\eta}\right) \tag{3.24}$$

The solution of (2.4), (2.5) satisfies the integral equation $u_2(t,z) = \exp(-\tilde{b} t)u_{2,0}(z) + \tilde{a}\int_0^t \exp(-\tilde{b}(t-\tau))u_1(\tau,z)d\tau$ for $(t,z) \in \mathfrak{R}_+ \times [0,1]$, from which we obtain the estimate for $t \ge 0$:

$$\|u_2[t]\|_{\infty,\eta} \le \exp(-\tilde{b} t)\|u_{2,0}\|_{\infty,\eta} + \tilde{b}^{-1}|\tilde{a}|\max_{0 \le s \le t}\left(\|u_1[s]\|_{\infty,\eta}\right) \tag{3.25}$$

Using Lemma 4.2 in [14] and (3.24), (3.25), we guarantee that there exists $\delta > 0$ (independent of the particular solution of (2.4), (2.5)) such that the following estimates hold for all $t \ge 0$:

$$\|u_1[t]\|_{\infty,\eta} \exp(\delta t) \le \|u_{1,0}\|_{\infty,\eta} + (1+\varepsilon)\frac{(1+\zeta^{-1})}{\sin(\theta)}\max_{0 \le s \le t}(|d(s)|\exp(\delta s))$$
$$+ \frac{(1+\varepsilon)|r|\tilde{b}(1+\zeta)}{K + (\pi - 2\theta)^2}\max_{0 \le s \le t}\left(\|u_2[s]\|_{\infty,\eta} \exp(\delta s)\right) \tag{3.26}$$

$$\|u_2[t]\|_{\infty,\eta} \exp(\delta t) \le \|u_{2,0}\|_{\infty,\eta} + (1+\varepsilon)\tilde{b}^{-1}|\tilde{a}|\max_{0 \le s \le t}\left(\|u_1[s]\|_{\infty,\eta} \exp(\delta s)\right) \tag{3.27}$$

Defining $p_i(t) := \max_{0 \le s \le t}\left(\|u_i[s]\|_{\infty,\eta} \exp(\delta s)\right)$ for $i = 1,2$ and $t \ge 0$, we obtain from (3.26), (3.27) the following inequalities for all $t \ge 0$:

$$p_1(t) \le \|u_{1,0}\|_{\infty,\eta} + \frac{(1+\varepsilon)|r|\tilde{b}(1+\zeta)}{K + (\pi - 2\theta)^2} p_2(t) + (1+\varepsilon)\frac{(1+\zeta^{-1})}{\sin(\theta)}\max_{0 \le s \le t}(|d(s)|\exp(\delta s)) \tag{3.28}$$

$$p_2(t) \le \|u_{2,0}\|_{\infty,\eta} + (1+\varepsilon)\tilde{b}^{-1}|\tilde{a}|p_1(t) \tag{3.29}$$

Inequalities (3.28), (3.29) in conjunction with (3.23) imply the following estimate for all $t \ge 0$:

$$p_1(t) + p_2(t) \leq (1-L)^{-1}\left(1+(1+\varepsilon)\tilde{b}^{-1}|\tilde{a}|\right)\|u_{1,0}\|_{\infty,\eta} + (1-L)^{-1}\left(1+\frac{(1+\varepsilon)|r|\tilde{b}(1+\zeta)}{K+(\pi-2\theta)^2}\right)\|u_{2,0}\|_{\infty,\eta}$$
$$+\left((1+\varepsilon)\tilde{b}^{-1}|\tilde{a}|+1\right)(1-L)^{-1}(1+\varepsilon)\frac{(1+\zeta^{-1})}{\sin(\theta)}\max_{0\leq s\leq t}\left(|d(s)|\exp(\delta s)\right)$$
(3.30)

Estimate (2.8) with appropriate constants $M, \gamma > 0$ is a direct consequence of (3.30), definitions $p_i(t) := \max_{0\leq s\leq t}\left(\|u_i[s]\|_{\infty,\eta}\exp(\delta s)\right)$ for $i = 1,2$ and the fact that there exist constants $0 < K_1 < K_2$ such that $K_1\|u\|_\infty \leq \|u\|_{\infty,\eta} \leq K_2\|u\|_\infty$ for all $u \in C^0([0,1])$. The proof is complete. ◁

**Proof of Theorem 2.4:** Without loss of generality we may assume that
$$|k| < 1 \tag{3.31}$$
If not then we can perform first the transformation $w_1(t,z) = u_1(t,z)$, $w_2(t,z) = u_2(t,z)/(1+|k|)$ and work with $w_1, w_2$ instead of $u_1, u_2$. It suffices to show that for every $T > 0$, $u_{1,0} \in D$, $u_{2,0} \in C^1([0,1])$ with $u_{2,0}(0) = ku_{1,0}(1)$, $\theta_{1,0} \in D$ $u'_{2,0}(0) = -c^{-1}k\theta_{1,0}(1)$, where $\theta_{1,0}(z) = pu''_{1,0}(z) + au_{1,0}(z) + \int_0^1 b(z,s)u_{2,0}(s)ds$ for $z \in [0,1]$, there exists a unique pair of mappings $u_1 \in C^0([0,T]\times[0,1])\cap C^1((0,T]\times[0,1])$, $u_2 \in C^1([0,T]\times[0,1])$ with $u_1[t] \in C^2([0,1])$ for $t \in (0,T]$ satisfying (2.20) and the equations

$$\frac{\partial u_1}{\partial t}(t,z) - p\frac{\partial^2 u_1}{\partial z^2}(t,z) - au_1(t,z) - \int_0^1 b(z,s)u_2(t,s)ds = \frac{\partial u_2}{\partial t}(t,z) + c\frac{\partial u_2}{\partial z}(t,z) = 0$$
$$\text{for } (t,z) \in (0,T]\times(0,1) \tag{3.32}$$

$$u_1(t,0) = \frac{\partial u_1}{\partial z}(t,1) - qu_1(t,1) = 0, \text{ for } t \in [0,T] \tag{3.33}$$

$$u_2(t,0) = ku_1(t,1), \text{ for } t \in [0,T] \tag{3.34}$$

Let $T > 0$, $u_{1,0} \in D$, $u_{2,0} \in C^1([0,1])$ with $u_{2,0}(0) = ku_{1,0}(1)$, $\theta_{1,0} \in D$ $u'_{2,0}(0) = -c^{-1}k\theta_{1,0}(1)$, where $\theta_{1,0}(z) = pu''_{1,0}(z) + au_{1,0}(z) + \int_0^1 b(z,s)u_{2,0}(s)ds$ for $z \in [0,1]$, be given (arbitrary).

It follows from Corollary 4.1 that any classical solution of (3.32), (3.33), (3.34), (2.20), $u_1 \in C^0([0,T]\times[0,1])\cap C^1((0,T]\times[0,1])$, $u_2 \in C^1([0,T]\times[0,1])$ with $u_1[t] \in C^2([0,1])$ for $t \in (0,T]$, must satisfy the following integral equations for all $(t,z) \in [0,T]\times[0,1]$:

$$u_1(t,z) = \sum_{n=1}^\infty \phi_n(z)\exp(-\lambda_n t)\int_0^1 \phi_n(s)u_{1,0}(s)ds + \sum_{n=1}^\infty \phi_n(z)\int_0^t \exp(-\lambda_n(t-\tau))\left(\int_0^1\int_0^1 \phi_n(s)b(s,l)u_2(\tau,l)dlds\right)d\tau$$
(3.35)

$$u_2(t,z) = \begin{cases} ku_1(t-c^{-1}z,1) \text{ for } ct > z \\ u_{2,0}(z-ct) \text{ for } ct \leq z \end{cases} \tag{3.36}$$

By virtue of (3.31), for sufficiently large $f > 0$ the mapping $\left(C^0([0,T]\times[0,1])\right)^2 \ni (v_1,v_2) \to (R_1(v_1,v_2), R_2(v_1,v_2)) \in \left(C^0([0,T]\times[0,1])\right)^2$ defined by

$$(R_1(v_1,v_2))(t,z) := \sum_{n=1}^\infty \phi_n(z)\exp(-(f+\lambda_n)t)\int_0^1 \phi_n(s)u_{1,0}(s)ds +$$
$$+ \sum_{n=1}^\infty \phi_n(z)\int_0^t \exp(-(f+\lambda_n)(t-\tau))\left(\int_0^1\int_0^1 \phi_n(s)b(s,l)v_2(\tau,l)dlds\right)d\tau$$
(3.37)

$$(R_2(v_1,v_2))(t,z) := \begin{cases} k\exp(-fc^{-1}z)v_1(t-c^{-1}z,1) \text{ for } ct > z \\ \exp(-ft)u_{2,0}(z-ct) \text{ for } ct \leq z \end{cases} \quad (3.38)$$

is a contraction. Therefore, there exists a unique $(v_1,v_2) \in (C^0([0,T]\times[0,1]))^2$ satisfying the equation $(v_1,v_2) = (R_1(v_1,v_2),R_2(v_1,v_2))$. By defining $u_i(t,z) = \exp(ft)v_i(t,z)$ for $i=1,2$, $(t,z)\in[0,T]\times[0,1]$, it follows from (3.37), (3.38) that there exists a unique $(u_1,u_2) \in (C^0([0,T]\times[0,1]))^2$ satisfying (3.35), (3.36) for $(t,z)\in[0,T]\times[0,1]$. Thus, the initial-boundary value problem (3.32), (3.33), (3.34), (2.20) has at most one classical solution.

Next, we construct the unique classical solution of the initial-boundary value problem (3.32), (3.33), (3.34), (2.20). Consider the integral equations for $(t,z)\in[0,T]\times[0,1]$:

$$\theta_1(t,z) = \sum_{n=1}^{\infty}\phi_n(z)\exp(-\lambda_n t)\int_0^1 \phi_n(s)\theta_{1,0}(s)ds + \sum_{n=1}^{\infty}\phi_n(z)\int_0^t \exp(-\lambda_n(t-\tau))\left(\int_0^1\int_0^1 \phi_n(s)b(s,l)\theta_2(\tau,l)dlds\right)d\tau \quad (3.39)$$

$$\theta_2(t,z) = \begin{cases} k\theta_1(t-c^{-1}z,1) \text{ for } ct > z \\ -cu'_{2,0}(z-ct) \text{ for } ct \leq z \end{cases} \quad (3.40)$$

Since $u_{2,0} \in C^1([0,1])$, $\theta_{1,0} \in D$ with $u'_{2,0}(0) = -c^{-1}k\theta_{1,0}(1)$ and by virtue of (3.31), it follows that for sufficiently large $f>0$, the mapping $(C^0([0,T]\times[0,1]))^2 \ni (q_1,q_2) \to (\tilde{R}_1(q_1,q_2),\tilde{R}_2(q_1,q_2)) \in (C^0([0,T]\times[0,1]))^2$ defined by

$$(\tilde{R}_1(q_1,q_2))(t,z) := \sum_{n=1}^{\infty}\phi_n(z)\exp(-(f+\lambda_n)t)\int_0^1 \phi_n(s)\theta_{1,0}(s)ds + \\ + \sum_{n=1}^{\infty}\phi_n(z)\int_0^t \exp(-(f+\lambda_n)(t-\tau))\left(\int_0^1\int_0^1 \phi_n(s)b(s,l)q_2(\tau,l)dlds\right)d\tau \quad (3.41)$$

$$(\tilde{R}_2(q_1,q_2))(t,z) := \begin{cases} k\exp(-fc^{-1}z)q_1(t-c^{-1}z,1) \text{ for } ct > z \\ -c\exp(-ft)u'_{2,0}(z-ct) \text{ for } ct \leq z \end{cases} \quad (3.42)$$

is a contraction. Therefore, there exists a unique $(q_1,q_2) \in (C^0([0,T]\times[0,1]))^2$ satisfying the equation $(q_1,q_2) = (\tilde{R}_1(q_1,q_2),\tilde{R}_2(q_1,q_2))$. By defining $\theta_i(t,z) = \exp(ft)q_i(t,z)$ for $i=1,2$, $(t,z)\in[0,T]\times[0,1]$, it follows from (3.41), (3.42) that there exists a unique $(\theta_1,\theta_2) \in (C^0([0,T]\times[0,1]))^2$ satisfying (3.39), (3.40) for $(t,z)\in[0,T]\times[0,1]$. It is a matter of straightforward calculations to verify (by using (3.39), (3.40) and the fact that $\int_0^1 \phi_n(s)(pu''_{1,0}(s) + au_{1,0}(s))ds = -\lambda_n\int_0^1 \phi_n(s)u_{1,0}(s)ds$) that the functions

$$u_i(t,z) = u_{i,0}(z) + \int_0^t \theta_i(\tau,z)d\tau, \quad i=1,2 \quad (3.43)$$

satisfy (3.35), (3.36) for $(t,z)\in[0,T]\times[0,1]$. Notice that (3.43) implies that $\frac{d}{dt}u_1(t,1) = \theta_1(t,1)$ is continuous on $[0,T]$. Consequently, definition (3.36) implies that $u_2 \in C^1([0,T]\times[0,1])$. Using Corollary 4.1 and (3.35) we obtain that $u_1 \in C^0([0,T]\times[0,1]) \cap C^1((0,T]\times[0,1])$ with $u_1[t] \in C^2([0,1])$ for $t \in (0,T]$. Moreover, (3.32), (3.33), (3.34), (2.20) hold. The proof is complete. ◁

**Proof of Lemma 2.5:** Consider the function $f(\omega) := \omega\cot(\omega)$ on $(0,\pi)$. It holds that $f'(\omega) = \frac{\cos(\omega)\sin(\omega) - \omega}{\sin^2(\omega)} = \frac{\sin(2\omega) - 2\omega}{2\sin^2(\omega)} < 0$ for all $\omega \in (0,\pi)$ and consequently $f(\omega)$ is strictly decreasing on $(0,\pi)$. Notice that $b_1 \in \left(-\frac{\pi}{2}, \frac{\pi}{2}\right)$ satisfies the equation $f\left(\frac{\pi}{2} - b_1\right) = q$. Therefore, for all $0 < \omega < \frac{\pi}{2} - b_1$ it holds that $\omega \cot(\omega) > q$. Moreover, since $p\left(\frac{\pi}{2} - b_1\right)^2 > a$ (recall (2.19)), the inequality $p\omega^2 > a$ may also be assumed to hold for $\omega < \frac{\pi}{2} - b_1$ sufficiently close to $\varphi = \frac{\pi}{2} - b_1$. Finally, since $\omega \cot(\omega) > q$, there exists $\theta \in (0, \pi - \omega)$ sufficiently small so that $\omega \cot(\omega + \theta) > q$. Consequently, the specifications $\theta \in (0, \pi)$, $\omega \in [0, \pi - \theta)$ with $p\omega^2 > a$ and $\omega \cot(\omega + \theta) > q$ hold. The proof is complete. ◁

We end this section by providing the proof of Theorem 2.6.

**Proof of Theorem 2.6:** By virtue of (2.21) there exists sufficiently small $\varepsilon > 0$ such that

$$(1+\varepsilon)^2 |k| \max_{0 \le z \le 1}\left(\frac{\sin(\omega+\theta)}{\sin(\omega z + \theta)} \int_0^1 |b(z,s)| ds\right) < p\omega^2 - a \quad (3.44)$$

Let $u_{1,0}, u_{2,0} \in C^1([0,1])$ with $u_{1,0}(1) = u'_{1,0}(1) - q u_{1,0}(1) = u_{2,0}(0) - k u_{1,0}(1) = 0$, for which there exists a unique pair of mappings $u_1 \in C^0(\Re_+ \times [0,1]) \cap C^1((0,+\infty) \times [0,1])$, $u_2 \in C^1(\Re_+ \times [0,1])$ with $u_1[t] \in C^2([0,1])$ for $t > 0$ satisfying (2.15), (2.16), (2.17), (2.20).

Define the positive function $\eta(z) := \sin(\theta + \omega z)$ for $z \in [0,1]$ and the norm $\|u\|_{\infty,\eta} := \max_{0 \le z \le 1}\left(\frac{|u(z)|}{\eta(z)}\right)$. Notice that Assumptions (H1), (H2), (H3), (H4) hold for the PDE problem $\frac{\partial u_1}{\partial t}(t,z) = p\frac{\partial^2 u_1}{\partial z^2}(t,z) + a u_1(t,z) + f(t,z)$ with $u_1(t,0) = \frac{\partial u_1}{\partial z}(t,1) - q u_1(t,1) = 0$. More specifically, Assumption (H4) holds with $\eta$ as defined above and $\sigma := p\omega^2 - a$. Applying Corollary 4.3 and recognizing that the solution of (2.15), (2.16) satisfies the equation $\frac{\partial u_1}{\partial t}(t,z) = p\frac{\partial^2 u_1}{\partial z^2}(t,z) + a u_1(t,z) + f(t,z)$ with $f(t,z) = \int_0^1 b(z,s) u_2(t,s) ds$, we obtain the following estimate for all $t \ge 0$:

$$\|u_1[t]\|_{\infty,\eta} \le \exp(-\sigma t) \|u_{1,0}\|_{\infty,\eta} + \sigma^{-1} B \max_{0 \le s \le t}(\|u_2[s]\|_\infty) \quad (3.45)$$

where $B := \max_{0 \le z \le 1}\left(\frac{1}{\eta(z)} \int_0^1 |b(z,s)| ds\right)$. The solution of (2.15), (2.17) also satisfies the equation (3.36) for $(t,z) \in \Re_+ \times [0,1]$, from which we obtain the estimate:

$$\|u_2[t]\|_\infty \le |k|\eta(1) \max_{0 \le s \le t}(\|u_1[s]\|_{\infty,\eta}) + M_\mu \exp(-\mu t) \|u_{2,0}\|_\infty, \text{ for all } t \ge 0 \text{ and } \mu > 0 \quad (3.46)$$

where $M_\mu := \exp(\mu c^{-1})$. Using Lemma 4.2 in [14] in conjunction with (3.45), (3.46) and the semigroup property, we guarantee that there exists $\delta \in (0,1)$ (independent of the particular solution of (2.15), (2.16), (2.17)) such that the following estimates hold for $t \ge 0$:

$$\|u_1[t]\|_{\infty,\eta} \exp(\delta t) \le \|u_{1,0}\|_{\infty,\eta} + (1+\varepsilon)\sigma^{-1} B \max_{0 \le s \le t}(\|u_2[s]\|_\infty \exp(\delta s)) \quad (3.47)$$

$$\|u_2[t]\|_\infty \exp(\delta t) \le (1+\varepsilon)|k|\eta(1) \max_{0 \le s \le t}(\|u_1[s]\|_{\infty,\eta} \exp(\delta s)) + \exp(c^{-1}) \|u_{2,0}\|_\infty \quad (3.48)$$

Defining $p_1(t) := \max_{0 \le s \le t}(\|u_1[s]\|_{\infty,\eta} \exp(\delta s))$, $p_2(t) := \max_{0 \le s \le t}(\|u_2[s]\|_\infty \exp(\delta s))$ for $t \ge 0$, we obtain from (3.47), (3.48) the following inequalities for all $t \ge 0$:

$$p_1(t) \le \|u_{1,0}\|_{\infty,\eta} + (1+\varepsilon)\sigma^{-1} B p_2(t), \quad p_2(t) \le (1+\varepsilon)|k|\eta(1) p_1(t) + \exp(c^{-1}) \|u_{2,0}\|_\infty \quad (3.49)$$

Inequalities (3.49) in conjunction with (3.44) (which together with definitions $\sigma := p\omega^2 - a$, $B := \max_{0 \leq z \leq 1}\left(\frac{1}{\eta(z)}\int_0^1 |b(z,s)|ds\right)$, $\eta(z) := \sin(\theta + \omega z)$ implies the inequality $(1+\varepsilon)^2 \sigma^{-1} B|k|\eta(1) < 1$) imply the following estimate for all $t \geq 0$:

$$p_1(t) + p_2(t) \leq \left(1 - (1+\varepsilon)^2 \sigma^{-1} B|k|\eta(1)\right)^{-1} \left((1+\varepsilon)\sigma^{-1} B + 1\right)\exp(c^{-1})\|u_{2,0}\|_\infty$$
$$+ \left(1 - (1+\varepsilon)^2 \sigma^{-1} B|k|\eta(1)\right)^{-1} \left(1 + (1+\varepsilon)|k|\eta(1)\right)\|u_{1,0}\|_{\infty,\eta} \qquad (3.50)$$

Estimate (2.22) with appropriate constant $M > 0$ is a direct consequence of (3.50), definitions $p_1(t) := \max_{0 \leq s \leq t}\left(\|u_1[s]\|_{\infty,\eta} \exp(\delta s)\right)$, $p_2(t) := \max_{0 \leq s \leq t}\left(\|u_2[s]\|_\infty \exp(\delta s)\right)$ and the fact that there exist constants $0 < K_1 < K_2$ such that $K_1 \|u\|_\infty \leq \|u\|_{\infty,\eta} \leq K_2 \|u\|_\infty$ for all $u \in C^0([0,1])$. The proof is complete. ◁

## 4. Concluding Remarks

We have shown that the small-gain methodology is useful for the stability analysis of hyperbolic-parabolic PDE loops. However, as it well-known from finite-dimensional systems, the small-gain analysis is usually a conservative methodology. It is not known whether the sufficient conditions provided by the main results of the present work are conservative or not, because we are not aware of any similar results in the spatial sup norm for the PDE loops that are studied in the present work. However, an extensive numerical study of the wave equation with Kelvin-Voigt damping and possible viscous damping with respect to sinusoidal boundary oscillations may reveal the level of conservatism for the estimation of the gain in Fig. 1. Future research work will address such issues.

The results of the present paper can be extended to the study of loops containing more than two PDEs. This is a topic for future research. It is also expected that the study of systems of PDEs with possible non-local reaction terms and boundary interconnections will open new research directions, because such systems may exhibit complicated dynamic behaviors.